\documentclass[reqno]{amsart}
\usepackage{amsfonts}
\usepackage{amssymb}
\usepackage{amsthm}
\usepackage{upref}
\usepackage{enumerate}
\makeatletter
\def\LaTeX{\leavevmode L\raise.42ex
    \hbox{\kern-.3em\size{\sf@size}{0pt}\selectfont A}\kern-.15em\TeX}
 \makeatother
 
\DeclareMathOperator{\clos}{clos}

\numberwithin{equation}{section}

\newtheorem{lemma}{Lemma}[section]
\newtheorem{theorem}[lemma]{Theorem} 
\newtheorem{corollary}[lemma]{Corollary}

\theoremstyle{definition}

\newtheorem{definition}[lemma]{Definition}

\newtheorem{assumption}[lemma]{Assumption}

\newtheorem{remark}[lemma]{Remark}

  \newcommand{\slim}{\operatorname{s-lim}}

  \newcommand{\e}{\eqref}

\newcommand{\q}{\quad}

\newcommand{\ov}{\overline}
\newcommand{\wt}{\widetilde}

\renewcommand\Im{\operatorname{Im}}
\renewcommand\Re{\operatorname{Re}}

\newenvironment{pf}{\begin{proof}}{\end{proof}}

\def\qqq{\mathrel{\subset\mkern-15mu\lower.38ex\hbox{${\scriptscriptstyle\rightarrow}$}}}

\let\cal\mathcal

\let\Bbb\mathbb

\begin{document}
\title 
{A note on the Schr\"odinger operator with a long-range potential}
\author{ D. R. Yafaev  }
\address{   Univ  Rennes, CNRS, IRMAR-UMR 6625, F-35000
    Rennes, France and SPGU, Univ. Nab. 7/9, Saint Petersburg, 199034 RUSSIA}
\email{yafaev@univ-rennes1.fr}
 

\thanks {Supported by  project   Russian Science Foundation   17-11-01126}

\begin{abstract}
Our goal  is to develop spectral and   scattering theories for the one-dimensional Schr\"odinger operator with a long-range potential $q(x)$, $x\geq 0$. Traditionally, this problem is studied with a help of the Green-Liouville approximation. This  requires conditions on the first two derivatives $q' (x)$ and $q'' (x)$. We suggest a new Ansatz that allows us to
develop a consistent theory under the only assumption $q' \in L^1$.  
   \end{abstract}

\maketitle

\thispagestyle{empty}

%***********************************************************
\section{Introduction}
%***********************************************************
 
\subsection{Short- and long-range potentials}

 The classical result of H.~Weyl  \cite{Weyl1, Weyl}  (see also the book \cite{Titch}) states that under very general circumstances  a differential  equation 
\begin{equation}
 -f''(x,z)+ q(x) f(x,z)= z f(x,z), \q  x\geq 0, \q \ov{q(x)}=q(x),  
 \label{eq:As1}\end{equation}
 where $ \Im z\neq 0 $, 
 has a solution $f(\cdot,z) \in L^2 ({\Bbb R}_{+})$. This fact, however, has no direct  spectral consequences for the 
 Schr\"odinger operator $H= -  d^2/dx^2 +q(x) $ (with some boundary condition at the point $x=0$)  in the space $L^2 ({\Bbb R}_{+})$,  except that its spectrum is simple. An advanced spectral analysis of the operator $H$ and scattering theory for the pair $H_{0} =-  d^2/dx^2$, $H$ require 
  the continuity of the solutions $f(\cdot,z)$ as $\Im z\to 0$ which can be verified only under
    some specific assumptions on the potential $q(x)$.

We suppose that  $q(x)\to 0$ as $x\to\infty$ and distinguish the short-  and long-range cases. In  the short-range case, it is assumed that  $q \in L^1 ({\Bbb R}_{+})$. This allows one to construct a solution    $\vartheta(x,z) $, known as the Jost solution,   of equation \e{eq:As1} with   asymptotics $\vartheta(x,z)\sim e^{-\sqrt{-z} x}$ as $x\to \infty$; it is  supposed here that $\Re \sqrt{-z}>0$ so that $\vartheta(\cdot,z) \in L^2 ({\Bbb R}_{+})$. It turns out that the function $\vartheta(x,z)$ is continuous in $z$ as $\Im z\to 0$. This fact is crucial for the  analysis of the operator $H$. It permits  (see, e.g., Sections 3.1 and 3.2 of the book \cite{Y}) to show that the   structure of its positive spectrum is essentially the same as that of the ``free" operator $H_{0}$. In particular, the positive spectrum of the operator $H$ is absolutely continuous. This fact follows from the continuity in an appropriate topology of the resolvent $R(z)=(H-z)^{-1}$ as $z$ approaches the positive spectrum of $H$. The last result is known as the limiting absorption principle. Note that
the Jost solutions are also widely discussed in the physics literature  (see, e.g., \S 1 and \S 2 in Chapter 12 of the book \cite{New}).

In the long-range case when $q \not\in L^1 ({\Bbb R}_{+})$, the definition of the Jost solution   has to be modified. It was shown by V.~B.~Matveev and M.~M.~Skriganov  in  \cite{MS} that, under some assumptions on the first two derivatives  $q' (x)$ and $q'' (x)$,   equation  \e{eq:As1} has a solution $\theta(x,z)$ described by the Green-Liouville Ansatz:
  \begin{equation}
\theta(x,z)\sim (q(x)-z)^{-1/4}\exp\Big( -\int_{0}^x (q(y)-z)^{1/2}dy\Big), \q x\to\infty;
 \label{eq:As2}\end{equation}
  it is  supposed here that $\Re\, (q(y)-z)^{1/2}>0$ so that  again  $\theta (\cdot,z) \in L^2 ({\Bbb R}_{+})$. We refer to the book \cite{Olver}, Chapter 6,  for a careful presentation of the   Green-Liouville  method. Given the existence of the solutions $\theta(x,z)$ and their continuity in $z$ up to the positive half-line, the spectral analysis of the operator $H$ is performed in  \cite{MS} essentially in the same way as in the short-range case.  The best possible conditions on $q(x)$ required by this method are  probably  $q' \in L^2$ and 
$q'' \in L^1$ (see \cite{Y-CMP}). 
  
  \subsection{Modified Green-Liouville Ansatz}
  
Our goal  is to develop spectral and stationary scattering theories for the Schr\"odinger operator $H$ with a long-range potential $q(x)$ under the only assumption $q' \in L^1 ({\Bbb R}_{+})$. Note that, for functions $q(x)$ not satisfying the short-range assumption $q\in L^1 ({\Bbb R}_{+})$, some conditions on  their derivatives are inavoidable. Indeed, the Wigner-von Neumann potential (see, e.g., Section~XIII.13 of the book \cite{RS}) has asymptotics $q(x)\sim x^{-1}\sin x$ as $x\to \infty$ and the corresponding operator $H$ has a positive eigenvalue. Thus, the positive spectrum of $H$ is not   absolutely continuous. More than that,
S.~N.~Naboko constructed in \cite{Nab} examples of Schr\"odinger operators with dense in $[0,\infty)$ point spectrum whose potentials decay only slightly slower than $x^{-1}$.

From analytic point of view, the present paper relies on a modification 
of the classical Green-Liouville Ansatz.  Actually, removing the factor  in front of the exponential, we replace the   formula \e{eq:As2} by a simpler one
\begin{equation}
\theta(x,z)\sim  \exp\Big( -\int_{0}^x (q(y)-z)^{1/2}dy\Big), \q x\to\infty.
 \label{eq:As3}\end{equation}
 The construction of    solutions $\theta(x,z)$ of equation \e{eq:As1} with such asymptotic behavior  under  the only assumption $q' \in L^1 ({\Bbb R}_{+})$ is the main new point of the paper. Then 
spectral and stationary scattering theories for the  operator $H$ can be developed along essentially the same lines as in the short-range case.

In the problem we consider, the Ansatz \e{eq:As3} is more efficient (and is much simpler)
than the Green-Liouville one.
However  the classical   Ansatz also has numerous advantages. For example, it was used in \cite{Y-CMP} to study low energy (as $z\to 0$) asymptotics of spectral and scattering data for potentials $q(x)$ decaying slower than  $x^{-2}$ as $x\to\infty$. Another important application of the Green-Liouville method is to differential equations \e{eq:As1} with  coefficients $q(x)$ tending to $+\infty$ or  to $-\infty$ as $x\to\infty$ (see  \cite{Olver}, Chapter 6).

  The approach we use works equally well for more general, than Schr\"odinger, differential operators. In the paper, we consider the operator
\begin{equation}
H= - \frac{d}{dx} p(x) \frac{d}{dx} +q(x)  , \q p(x)>0, \q \ov{q(x)}= q(x),
\label{eq:A}\end{equation}
with some boundary condition at the point $x=0$ in the space $L^2 ({\Bbb R}_{+})$. We choose the condition $f(0)=0$.

\subsection{Other methods}

The limiting absorption principle under   assumptions very close to $q'\in L^1$ was obtained long ago by the powerful Mourre method \cite{Mo1}; we refer to \S 6.9 of \cite{Y} where the conditions on $q$ were stated explicitly. The Mourre method works for very general operators (for example, for the Schr\"odinger operator in all dimensions), but the only spectral information it gives is the absence of the singular continuous spectrum. 
For example, it says nothing about multiplicity of the spectrum.

There are also specific one-dimensional methods   for a proof of the absolute continuity of the positive spectrum.
Thus  J.~Weidmann in \cite{Weidmann}, see Theorem~14.25,  proved this fact 
by the Gilbert-Pearson method \cite{GD}  under somewhat more general assumptions on   $q(x)$ than $q'\in L^1$.
 
The analytic approach we use here allows us to perform spectral analysis of the operator $H$   in a much more detailed way.

  \subsection{Structure of the paper}

Section~2 is central. Here we construct modified Jost solutions, introduce a multiplicative change of variables and reduce   differential equation \e{eq:As1} to a Volterra integral equation.

We study   a regular solution $\varphi(x,z)$ of equation \e{eq:As1} in Section~3. In particular,  we find its asymptotic behavior    as $x\to\infty$. Of course, the answers are quite different for $z>0$ and $z\not\in [0,\infty)$. The result for  $z\not\in [0,\infty)$ seems not to be   well known even in the short-range case.

Given the results of Sections~2 and 3, we follow the standard approach to   spectral analysis of the operator $H$   in Section~4. First, we define $H$  as a self-adjoint operator. Then we   build its eigenfunctions of the continuous spectrum and establish an expansion theorem over these eigenfunctions. The limiting absorption principle is a  by-product of these considerations.

Here is the list of miscellaneous results of Section~5: inclusion of an additional short-range term
$q_{\rm sr}\in L^1$, a general boundary condition at the point $x=0$, the problem on the whole line.

 Occasionally, the dependence of various  functions on $x$ and $z$ is   omitted in notation; $c$ and $C$ are different constants whose precise values are of no importance. We also use notation: $(A)_{+}= (|A|+A)/2$ for $A\in {\Bbb R}$.

%***********************************************************
\section{Modified Jost solutions}
%***********************************************************

\subsection{ Ansatz}

Let us consider a more general than \e{eq:As1} differential equation
\begin{equation}
 - ( p (x) f '(x, z ) )'+q(x) f (x ,z) = z  f (x ,z)
\label{eq:A1}\end{equation}
for $z\in {\Bbb C }\setminus {\Bbb R}=:\Pi$. We admit also that the spectral parameter $z$ belongs to the closure $\clos\Pi$ of $\Pi$, that is, $z=\lambda\pm i0$ where $\lambda\in {\Bbb R }\setminus  \{0\}$.
With respect to the  functions  $p(x)>0$ and $q(x) =\ov{q(x)}$, we accept the following

\begin{assumption}\label{PQ}
 \begin{enumerate}[{\rm(i)}]
\item
 The function $p(x)$ is absolutely continuous on  $ {\Bbb R}_{+}$.  The function $q(x)$  is absolutely continuous for $x\geq x_{0}$ where  $x_{0}$ may be arbitrary  large and   $q \in L^1(0,x_{0})$.
 \item
 The derivatives
\[
 p ' \in L^1 ({\Bbb R}_+),\quad q' \in L^1(x_{0},\infty)
\]
\item
The limits
 \begin{equation}
  \lim_{x\to\infty} p (x )=: p_{0}  >0, \q \lim_{x\to\infty} q (x ) =0 .
\label{eq:A2a}\end{equation}
 \end{enumerate}
\end{assumption}
 
It is sufficient to construct solutions of the  differential equation
 \e{eq:A1}  for large $x$. Then they can be standardly extended to all $x\geq 0$. A more general situation of this type     is discussed in Subs.~5.1.
 
Our goal in this section is to distinguish a solution $\theta (x ,z) $ of   equation \e{eq:A1} by its behavior as $x\to\infty$.
 To define it, we first
 exhibit  an explicit   function $a (x ,z) $ such that the remainder 
 \begin{equation}
r(x ,z): =- a (x ,z)^{-1}( p(x) a '(x, z ) )'+ q(x)-z  
\label{eq:A3}\end{equation}
in  equation \e{eq:A1}
is in $L^1$ (at infinity). Let us seek  $a (x ,z) $ in the form
 \begin{equation}
 a (x ,z)=e^{-\Omega (x ,z)}
\label{eq:A4}\end{equation}
and put
$
\omega (x ,z) =  \Omega' (x ,z)  .
$
Since
\begin{equation}
a(x ,z)^{-1} a '(x, z ) =-\omega (x ,z) ,
\label{eq:A6a}\end{equation}
we can rewrite \e{eq:A3} as
 \begin{align}
r(x ,z) &=  a (x ,z)^{-1}( p(x) \omega (x, z ) a (x, z ) )'+ q(x)-z  
\nonumber\\ 
&=  (p(x)  \omega(x ,z))'  -p(x)\omega(x ,z)^2+  q(x)-z .
\label{eq:A5}\end{align}

Let  us set
 \begin{equation}
\omega (x ,z) =  \sqrt{\frac{q(x)-z  } {p(x)} },
\label{eq:A7}\end{equation}
where we suppose that $\Re \omega (x ,z)>0$ for all $z\in \Pi$. 
Clearly,  the function $\omega (x ,z)$ is analytic in $z\in \Pi$ and continuous up to the cut along $\Bbb R$.
For $z\neq 0$, we choose $x_{1}= x_{1}(z)$ such that $|q(x)|\leq |z| /2$  for all $x\geq x_{1}$. 
Then  estimates
 \begin{equation}
0 < c\leq |\omega (x ,z) |\leq C<\infty
\label{eq:A7a}\end{equation}
and
 \begin{equation}
|\omega'(x ,z) |\leq C (| p'(x) | +|q'(x)|)
\label{eq:A7b}\end{equation}
are true
for all $x\geq x_{1}$, and hence $\omega' (\cdot,z)\in L^ 1 (x_{1},\infty)$ under Assumption~\ref{PQ}.  Here and below all estimates are uniform in $z$ from compact subsets of $\clos\Pi\setminus \{0\}$ (including the values of $z=\lambda\pm i0$ on the cut).

For the choice \e{eq:A7}, the remainder \e{eq:A5} equals
 \begin{equation}
r(x ,z) = (p(x)  \omega(x ,z))'  .
\label{eq:A8}\end{equation}
In view \e{eq:A7a}, \e{eq:A7b}, this yields the following result.

 \begin{lemma}\label{eik}
 Define the functions $\Omega (x ,z)$ and $a (x ,z)$ by the formulas
 \begin{equation}
 \Omega (x ,z)= \int_{0}^x  \sqrt{\frac{q(y)-z  } {p(y)} } dy
\label{eq:A11}\end{equation}
and \e{eq:A4}. Then the remainder \e{eq:A3} is given by formula \e{eq:A8} and  
$
r\in L^1 (  x_{1} (z),\infty).
$
 \end{lemma}

We emphasize that  the classical Green-Liouville Ansatz differs from \e{eq:A4} by the additional factor $\omega(x ,z)^{-1/2}$ in the right-hand side.

Since $\Re \omega (x,z)\geq  0$, it follows from \e{eq:A4}, \e{eq:A11} that 
\begin{equation}
\Big|  \frac{a (y,  z)} {a (x,  z)}\Big| \leq 1, \q y\geq x.
\label{eq:A16x}\end{equation}
Moreover,  $\Re \omega (x,z)\geq c(z)>0$ for $\Im z\neq 0$, so that we have a stronger estimate
\begin{equation}
\Big|  \frac{a (y,  z)} {a (x,  z)}\Big| \leq e^{-c(z) (y-x)}, \q y\geq x, \q c(z)>0, \q \Im z\neq 0.
\label{eq:A16X}\end{equation}

\subsection{Multiplicative substitution}
Instead of a solution $ \theta (x ,z)$ of equation \e{eq:A1}, we introduce a  function
\begin{equation}
 u (x ,z)= a(x ,z)^{-1} \theta (x ,z).
\label{eq:A10}\end{equation}

\begin{lemma}\label{subst}
Let $r (x ,z) $ be given by formula \e{eq:A8}. Then
 equation  \e{eq:A1} for $ \theta(x ,z)$ is equivalent to the equation
\begin{equation}
 ( p (x) u '(x, z ) )'- 2 \omega (x ,z) p (x) u '(x, z ) = r (x ,z) u (x, z )
\label{eq:A12}\end{equation}
for the function \e{eq:A10}. 
 \end{lemma}

\begin{pf}
A direct differentiation of the relation $\theta = a u$ shows that
\[
(p\theta' )'=   (p a' u)'+ (p a u')'
= (p a' )' u + (p u' )' a -2 \omega pa u'
\]
where we have taken \e{eq:A6a} into account. It follows that
\[
-(p\theta' )' + (q-z)\theta= \big(-(pa' )' + (q-z)a  \big) u + \big(-(pu' )' + 2 \omega p u'  \big) a.
\]
Since the first term on the right equals $ rau$, we see  that
\[
-(p\theta' )' + (q-z)\theta= \big(-(pa' )' + (q-z)a  \big) u + \big(-(pu' )' + 2 \omega p u'  + ru\big) a.
\]
 Thus equations \e{eq:A1}  and \e{eq:A12} coincide.
\end{pf}

Next, we reduce differential equation \e{eq:A12} to an integral equation. 

\begin{lemma}\label{int}
Let $z\in\Pi$, and let $u (x, z )$ be a solution of  differential equation \e{eq:A12} such that
\begin{equation}
 \lim_{x\to\infty}u (x, z ) =1.
\label{eq:A12a}\end{equation}
Then
\begin{equation}
 u(x, z )= 1 + \int_{x}^\infty G (x,y,z)  r (y, z )u (y, z ) dy  
\label{eq:A17}\end{equation}
where 
\begin{equation}
   G (x,y,z)=  a (y,z)^2   \int_{x}^y p(s)^{-1} a (s,z)^{-2} ds.
\label{eq:A18}\end{equation}
 \end{lemma}

\begin{pf}
 Set
\begin{equation}
 v(x, z )=p(x) u '(x, z ) \q {\rm and} \q \rho (x, z ) = r (x ,z) u (x, z ).
\label{eq:A13}\end{equation}
Then  \e{eq:A12} yields a differential equation 
\[
  v '(x, z ) - 2 \omega (x ,z)   v(x, z ) = \rho(x ,z)  
\]
of first order for the function $v(x, z )$. In view of \e{eq:A6a},  its solution is given by the equality
\[
 v(x, z )= -  a(x,z)^{-2} \Big(\int _{x} ^\infty a(y,z)^2 \rho (y, z ) dy +c\Big)
\]
where we have to choose $c=0$ because $a(x,z)^{-2}$ exponentially grows as $x\to\infty$.
Therefore the function $u (x, z )$ satisfying \e{eq:A12a} and \e{eq:A13} can be recovered by the formula
\[
 u(x, z )=1 -\int_{x}^\infty p(s)^{-1}  v(s, z ) ds=1 + \int_{x}^\infty p(s)^{-1} a (s,z)^{-2} \Big(\int _{s} ^\infty a (y,z)^2 \rho (y, z ) dy \Big) ds.
\]
By virtue of \e{eq:A16X}, Fubini's theorem allows us to interchange    the order of integrations here.
This yields equation 
\e{eq:A17}.
\end{pf}

\subsection{Integral equation}
The following assertion plays the crucial role in the analysis of equation \e{eq:A17} as $z$ approaches the half-axis $(0,\infty)$ (the continuous spectrum of $H$).

 \begin{lemma}\label{eik1}
For  $z\in\clos\Pi\setminus\{0\}$ and  $y\geq x\geq x_{1} (z)$, kernel \e{eq:A18} is uniformly bounded:
 \begin{equation}
| G (x,y,z) |\leq C<\infty.
\label{eq:A19}\end{equation} 
 \end{lemma}

 \begin{pf}
 Set
 $ \tau  = (p \omega )^{-1}$.  According to \e{eq:A7a} and \e{eq:A7b}, we have
 \begin{equation}
\tau \in L^\infty \q {\rm and} \q \tau ' \in L^1
\label{eq:A18T}\end{equation}
   Integrating by parts, we see that 
 \[
  2     \int_{x}^y p(s)^{-1}  a (s)^{-2} ds=  \int_{x}^y   \tau (s )   d a (s)^{-2}
  =
\tau  (y)  a (y)^{-2} -  \tau  (x)  a (x)^{-2} -\int_{x}^y   \tau '(s)   a (s)^{-2}ds.
\]
Multiplying this equality by $a (y)^{2} $ and using  estimate \e{eq:A16x} and relations \e{eq:A18T}, we get bound \e{eq:A19}.
  \end{pf}

  Lemmas~\ref{eik} and \ref{eik1} allow us to solve the Volterra equation \e{eq:A17} by iterations. Let us state the corresponding result.
  
  \begin{lemma}\label{eik2}
  For  $z\in\clos\Pi\setminus\{0\}$,
equation  \e{eq:A17} has a $($unique$)$ bounded solution $u(x, z )$. This function obeys an estimate
 \begin{equation}
| u(x, z )-1|\leq C \varepsilon(x)
\label{eq:A20}\end{equation} 
where
 \[
\varepsilon(x)= \int_{x}^\infty \big(| p'(y) | +|q'(y)|\big) dy.
 \]
  \end{lemma}
 
 The next assertion is converse to Lemma~\ref{int}.
 
  \begin{lemma}\label{int1}
  For  $z\in\clos\Pi\setminus\{0\}$, a solution $u (x,z )$ of   integral  equation  \e{eq:A17} satisfies also differential  equation \e{eq:A12}.
 \end{lemma}
 
  \begin{pf}
 According to  \e{eq:A18} we have
\[
   G (x,x)= 0 \q {\rm and} \q  G_{x}' (x,y)= - p(x)^{-1}a (x  )^{-2}a (y )^2.
\]
Therefore it follows from \e{eq:A17}    that
\begin{equation}
 u'(x  )= - p(x)^{-1}a (x  )^{-2}  \int_{x}^\infty a (y  )^2  r (y )u (y  ) dy  
\label{eq:A17a}\end{equation}
and hence, by \e{eq:A6a},
\begin{equation}
(p(x) u'(x ))'=  -2a (x  )^{-2}  \omega (x  )\int_{x}^\infty a(y)^2  r(y  )u (y  ) dy  +
r(x  )u (x  ).
\label{eq:A17b}\end{equation}
Substituting expressions \e{eq:A17a} and \e{eq:A17b} into the left-hand side of \e{eq:A12}, we see that it equals
  $r (x  )u (x  )$. 
   \end{pf}
  
  Now  we are in a position to give a precise
  
  \begin{definition}\label{MJS}
    For  $z\in\clos\Pi\setminus\{0\}$,
  define the function $a (x, z )$ by formulas \e{eq:A4} and \e{eq:A11}. 
  Denote by $u(x, z )$   the function constructed in Lemma~\ref{eik2}. 
  The  (modified) Jost solution of equation \e{eq:A1} is defined by the formula
  \begin{equation}
   \theta(x, z )= a(x, z ) u(x, z )
   \label{eq:Jost}\end{equation} 
 for $x\geq x_{1} (z)$,  and then $   \theta(x, z )$ is extended to all $x\geq 0$ as a solution of equation \e{eq:A1}.
 \end{definition}
 
 It follows from \e{eq:A20} that
  \begin{equation}
 \theta(x, z )=a(x, z )\big(1 + O( \varepsilon(x))\big), \q x\to\infty.
\label{eq:A22}\end{equation} 
 Note also that
  \[
 \theta(x, \bar{z} )=\ov{\theta(x, z )}
 \]
 and, in particular, 
  \begin{equation}
 \theta(x, \lambda- i0 )=\ov{\theta(x, \lambda+ i0 )}, \q \lambda\in {\Bbb R}\setminus\{0\}.
 \label{eq:Acx}\end{equation}

  Let us summarize the results obtained.
 
  \begin{theorem}\label{EIK}
   Let Assumption~\ref{PQ} be satisfied, and let $z\in\clos\Pi\setminus\{0\}$.
Denote by $u(x, z )$   the function constructed in Lemma~\ref{eik2}. Then  the function $\theta(x, z )$ defined by equality \e{eq:A10} satisfies equation \e{eq:A1},  and it has  asymptotics
\e{eq:A22}. For every $x\geq 0$, the function $\theta(x, z )$ is analytic in $z\in  \Pi$  and is continuous up to the cut along $\Bbb R$
with possible exception of the point $z=0$. Asymptotics \e{eq:A22} is  uniform in $z$ from compact subsets    of the set $  \clos\Pi\setminus\{0\} $.  
 \end{theorem}
 
   \begin{corollary}\label{EIKc}
   For $\lambda>0$, set  
    \begin{equation}
 \Phi (x,\lambda)= \int_{0}^x \sqrt{\Big( \frac{\lambda-q(y)}{p(y)}\Big)_{+}}dy, \q
  K(\lambda)= \exp\Big(-\int_{0}^\infty \sqrt{\Big( \frac{ q(y)-\lambda}{p(y)}\Big)_{+}}dy\Big).
 \label{eq:K}\end{equation}
 Then
     \[
 \theta (x,\lambda\pm i0)=  
  K(\lambda)\exp\Big(\pm i \Phi (x,\lambda)\Big) (1+\varepsilon (x))\q\mbox{as}\q x\to\infty.
   \]
 \end{corollary}
 
 Let us make several additional observations.
 
  \begin{remark}\label{EIKrem}
  \begin{enumerate}[{\rm(i)}]
 \item
 Unlike the Jost solution in the short-range case, the function $\theta(x, z )$ is not analytic in the whole half-plane $\Re z< 0$ because, in general, $ \theta(x, \lambda- i0 )\neq \theta(x, \lambda+ i0 )$.  This circumstance is, however, inessential. In particular, it follows from \e{eq:Acx} that $ \theta(x, \lambda- i0 )=0$ if and only if $ \theta(x, \lambda+ i0 )=0$. This subject is further discussed in Subsect.~4.2.
  \item
  For $z\in\Pi$,
 relation \e{eq:A22}  distinguishes a 
unique solution of equation \e{eq:A1}.   Indeed,  the differential equations \e{eq:A1} and the integral equation \e{eq:A17} are equivalent,  and Lemma~\ref{eik2} ensures that the solution of \e{eq:A17} satisfying  \e{eq:A20} is unique.
   \item
    According to \e{eq:A17a} $u' (x)= O (\varepsilon(x))$ as $x\to \infty$, and hence the derivative $\theta'= a (-\omega u+ u')$ has asymptotics 
    \begin{equation}
 \theta'(x, z )= - \sqrt{-z/p_{0}}\, a(x, z )\big(1 + O( \varepsilon(x))\big), \q x\to\infty.
\label{eq:A22D}\end{equation}
\item
It follows from definition \e{eq:A7} that
 \[
 \omega (x ,z)=     \sqrt{ -z  / p_0 } +o(1) \q \mbox{   whence } \q \Omega (x ,z)=  x  \sqrt{ -z  / p_0 } +o(x)
\]
as $x\to\infty$. In particular, 
 \[
\theta (\cdot , z )\in L^2 ({\Bbb R}_{+})  \q \mbox{  for } \q z\in\Pi.
\]
  \end{enumerate}
 \end{remark}

   \subsection{Non-uniqueness of Jost solutions}

For long-range perturbations, there is no canonical choice of the Jost solution $\theta(x,z)$: one can replace 
$\theta(x,z)$ by 
\begin{equation}
\wt{\theta}(x,z)=\theta(x,z) b(z),  \q  b(\bar{z})=\ov{b(z)},
\label{eq:JM}\end{equation}
 where  $b(z)$  is some function analytic in $\Pi$ and continuous up to the cut along $\Bbb R$. 
 Thus the function $\Omega (x,z)$ can be replaced in \e{eq:A4}  by a function $\wt{\Omega}  (x,z)$ provided the difference $\Omega (x,z)-\wt{\Omega}  (x,z)$ has a finite limit as $x\to\infty$.
 
  This observation allows one to simplify expression \e{eq:A11}   if additional information on decay of $p_{1}(x)=p (x)-p_0$ and of $q(x)$ is available. Note, first,   that in the short-range case when $p_{1}\in L^1$ and $q \in L^1$ one can choose
$\Omega_{0} (x,z)= x \sqrt{-z/p_{0}}$ instead of $\Omega  (x,z)$. Indeed, in this case we have
\[
\lim_{x\to\infty}\big( \Omega  (x,z)-\Omega_{0} (x,z)\big)=\int_{0}^\infty\frac{p_{0} q(x) +z p_{1} (x)}
{p_{0}\sqrt{(q(x)-z)p(x)}  + p(x)\sqrt{-zp_{0}}} dx=:\beta(z).
\]
Recall   that in the short-range case the standard (non-modified) Jost solution $\vartheta(x,z)$ is distinguished by the asymptotics $\vartheta (x,z)\sim e^{- x \sqrt{-z/p_{0}}}$ as $x\to \infty$. Therefore we have  
$\theta(x,z)=e^{-\beta(z)}\vartheta (x,z)$.

 If  $p_{1}\in L^2$ and $q \in L^2$ (but $p_{1}\not\in L^1$, $q\not \in L^1$), we set
\begin{equation}
\Omega_1 (x,z)= \sqrt{  -z/p_0}\Big( x- (2z)^{-1} \int_{0}^x q(y) dy- (2p_{0})^{-1} \int_{0}^x p_{1}(y) dy \Big).
\label{eq:Ras2}\end{equation}
Then $\Omega  (x,z)-\Omega_1 (x,z)$ has  a finite limit as $x\to\infty$ so that function \e{eq:A11} can be replaced by \e{eq:Ras2}. 
Similarly, in  the case  $p_{1}\in L^3$ and $q \in L^3$ (but $p_{1}\not\in L^2$, $q\not \in L^2$),
  we set
\begin{equation}
\Omega_2 (x,z)=\Omega_1 (x,z)+  8^{-1}\sqrt{  - z /p_0}\Big(z^{-2} \int_{0}^x q(y)^2 dy +  3 p_{0}^{-2} \int_{0}^x p_{1}(y)^2 dy \Big)
\label{eq:Ras3}\end{equation}
where $\Omega_1 (x,z)$ is given by \e{eq:Ras2}. Then $\Omega  (x,z)-\Omega_2 (x,z)$ has  a finite limit as $x\to\infty$ so that function  \e{eq:A11} can be replaced by \e{eq:Ras3}. This procedure can be continued to treat a general case where  $p_{1}\in L^n$ and $q \in L^n$  for some integer $n$.

  \section{Regular solution}
  
   \subsection{Asymptotics at infinity}

 In addition to $\theta(x,z)$, we introduce a regular    solution $\varphi(x,z)$  of equation \e{eq:A1} distinguished by   conditions at the point $x=0$:
\begin{equation}
\varphi(0,z)=0, \q \varphi'(0,z)=1.
\label{eq:B}\end{equation}
For every $x\geq 0$, the function $\varphi(x,z)$ is analytic in $z\in \Bbb C$.

 For $z\in \clos\Pi \setminus\{0\}$, let us consider the Wronskian
\begin{equation}
w(z)= \{\varphi(\cdot,z) ,\theta(\cdot,z) \}:=p(x) (\varphi'(x,z) \theta(x,z) -\varphi(x,z) \theta'(x,z) )
\label{eq:Bw}\end{equation}
of the regular $\varphi(x,z)$ and Jost $ \theta(x,z)$ solutions of equation \e{eq:A1}. Since the right-hand side of \e{eq:Bw} does not depend on $x \geq 0$,  we can set $x= 0$ whence
\begin{equation}
w(z) =p(0)   \theta(0,z) .
\label{eq:Bw1}\end{equation}
  
Let us now consider the Jost solutions $\theta(x, \lambda\pm i0)$  on the cut along $(0,\infty)$.  
Using    Corollary~\ref{EIKc}  and calculating the
  Wronskian of $\theta(x,\lambda+i0)$ and $ \theta (x,\lambda-i0)$ for $x\to\infty$, we find that
\[
 w_{0} (\lambda) : =  \{ \theta(\cdot,\lambda+i0) , \theta(\cdot,\lambda- i0)\}  = 2i\sqrt{p_{0}\lambda}K(\lambda)^2.
\]
In particular, the Jost solutions are linearly independent.
It is easy to see that
\begin{equation}
\varphi(x,\lambda)=  \frac{\theta(x,\lambda+i0)w(\lambda-i0)-\theta(x,\lambda-i0)w(\lambda+i0)}{ \{ \theta(\cdot,\lambda+i0) , \theta(\cdot,\lambda- i0)\} }.
\label{eq:R1}\end{equation} 
Indeed, by \e{eq:Bw1}, the right-hand side of \e{eq:R1} equals $0$   for $x=0$ and    its derivative in $x$  equals 
$ 1$  for $x=0$. Thus the right-hand side of \e{eq:R1} satisfies equation \e{eq:A1} where $z=\lambda$ and conditions \e{eq:B}.
It follows from \e{eq:R1} that 
\[
w(\lambda\pm i0)\neq 0, \q \lambda>0,
\]
 since otherwise we would have $\varphi(x,\lambda)=0$ for all $x\geq 0$.

Now we set
\begin{equation}
\kappa(\lambda)=|w(\lambda\pm i0)|\q {\rm and}\q w(\lambda + i0)=\kappa(\lambda)e^{i\eta(\lambda)}.
\label{eq:J}\end{equation} 
By analogy with the short-range case, we use the terms ``limit amplitude" for $\kappa(\lambda)$ and ``limit phase" 
(or  scattering phase, or     phase shift) for $\eta(\lambda)$.
 According to Theorem~\ref{EIK} the amplitude $\kappa(\lambda)$ is a continuous function of  $\lambda>0$. The phase $\eta (\lambda)$ is defined by equations \e{eq:J} up to a term $2\pi n$ where $n$ is an integer. Since $\kappa(\lambda)> 0$, the  function  $\eta (\lambda)$ can also be chosen continuous in $\lambda>0$.

The next results is a direct consequence of Corollary~\ref{EIKc} and representation \e{eq:R1}.

\begin{theorem}\label{REG}
 Let Assumption~\ref{PQ} be satisfied.
For   each $\lambda>  0$, the regular soution $\varphi(x,\lambda)$ of equation \e{eq:A1} has asymptotics
\begin{equation}
\varphi(x,\lambda) =\frac{ \kappa(\lambda)}{ \sqrt{p_{0}\lambda}\,  K(\lambda)}\sin\big( \Phi(x,\lambda )-\eta (\lambda)\big)+\varepsilon(x)
\label{eq:Ras}\end{equation} 
 as $x\to\infty$ 
In particular, the operator $H$ has no positive eigenvalues. 
 \end{theorem}
 
 \begin{remark}\label{nu}
 In the short-range case, the Jost solution $\vartheta(x,\lambda\pm i0)$ is distinguished by the asymptotics $\vartheta(x,\lambda\pm i0)\sim e^{\pm ix\sqrt{\lambda/p_0}}$ as $x\to\infty$. Then the limit amplitude $\kappa(\lambda)$ and the limit phase $\eta(\lambda)$ are defined by relations \e{eq:J} where $w(\lambda)=\{\varphi (\cdot, \lambda), \vartheta  (\cdot, \lambda+ i0)\}$.
 \end{remark}
 
  \begin{remark}\label{nu1}
Expressions \e{eq:R1} and hence \e{eq:Ras} are of course   invariant with respect to the change \e{eq:JM} of the Jost solution. This means that the amplitude factor and the phase $\Phi(x,\lambda )-\eta (\lambda)$ in \e{eq:Ras} do not depend on the regularization of the Jost solution. However we emphasize that, separately,  the terms $\Phi(x,\lambda ) $ and $ \eta (\lambda)$ do depend on it.
  Hence, in the long-range case, the definition of the scattering phase $\eta (\lambda)$ is not intrinsic. 
 \end{remark}

 \subsection{Exponentially growing solutions}
 Generically,  for $z\not\in [0,\infty)$,  the regular solution $\varphi(x,z)$ grows exponentially as $x\to\infty$.
 Here we find its asymptotics. The method below was used in the short-range case in  \S 4.1 of the book \cite{Y} but seems not to be commonly known.

Let     $\theta (x,z)$ be the Jost solution of equation  \e{eq:A1}.  For $z=\lambda<0$, we can pick either $\theta (x,\lambda+i0)$ or $\theta (x,\lambda-i0)$.
 We choose $\varrho=\varrho (z)$   in such a way that $\theta(x,z) \neq 0$ for all $x\geq \varrho (z)$. If $\Im z\neq 0$, we can set $\varrho (z)=0$ because the equality
$\theta(x_{0},z) = 0$ would imply that the self-adjoint operator \e{eq:A} in the space $L^2
(x_{0},\infty)$ with  the boundary condition $f(x_{0})=0$ has   complex eigenvalue $z$. 
Let us   introduce an exponentially growing solution $\xi (x,z)$ of  equation  \e{eq:A1}.

\begin{theorem}\label{EG}
Let Assumption~\ref{PQ} be satisfied, and
let $z\in {\Bbb C}\setminus [0,\infty)$. Then the function
  \begin{equation}
 \xi (x,z)= \theta (x,z) \int_{\varrho (z)}^x  \theta(y,z)^{-2} p(y)^{-1} dy , \q x\geq \varrho (z),
\label{eq:D}\end{equation}
   satisfies equation  \e{eq:A1} and
 \begin{equation}
 \xi (x,z)= \Big( 2  \sqrt{-p_{0}z} \: a (x,z) \Big)^{-1} (1+ o(1)), \q  \xi' (x,z)= \Big( 2  p_{0}   a (x,z) \Big)^{-1} (1+
 o(1))
\label{eq:D4}\end{equation}
as $x\to\infty$.
 \end{theorem}

\begin{pf}
Differentiating \e{eq:D}, we find that
\[
-(p (x)\xi' (x))'+ (q(x)-z)\xi (x)= \Big( -(p (x) \theta'(x))'+ (q(x)-z)\theta (x)\Big)\int_{\varrho }^x  \theta (y)^{-2} p(y) ^{-1} dy 
\]
which implies equation  \e{eq:A1} for $\xi (x)$.

Integrating by parts, we see that
 \begin{equation}
2  \int_{\varrho}^x \theta(y )^{-2} p(y)^{-1}dy
  =
 t  (x) a (x)^{-2} -  t  (\varrho) a (\varrho)^{-2}
   -   \int_{\varrho}^x   t '(y)   a (y)^{-2}d y
\label{eq:A18xN}\end{equation}
where $ t  = (p \omega u^2)^{-1}$. Let us multiply this equality by $a(x) \theta(x)$ and
consider the limit $x\to\infty$. It follows from Lemma~\ref{eik2} and Theorem~\ref{EIK} that
the first term on the right 
\begin{equation}
 t  (x)  \theta(x) a (x)^{-1}= (-p_{0} z )^{-1/2} + O(\varepsilon(x)).
\label{eq:A8xN}\end{equation}
  The second term   tends to $0$ exponentially. Finally, using estimate \e{eq:A16X}, we find that
 \begin{equation}
\Big| a(x) \theta(x)  \int_{\varrho}^x   t '(y)   a (y)^{-2}d y\Big| \leq C\int_{0}^x e^{-c(z) (x-y)} |t '(y) | dy.
\label{eq:A8xM}\end{equation}
    Since according to Lemma~\ref{eik2} (see also equality \e{eq:A17a}) $t' (x)\to 0$ as $x\to \infty$, the same is true for expression    \e{eq:A8xM}.  
Therefore relations \e{eq:A18xN} and \e{eq:A8xN} imply asymptotic formulas \e{eq:D4}.
   \end{pf}

 Using asymptotics \e{eq:A22}, \e{eq:A22D} and \e{eq:D4},
 we can calculate the Wronskian of the solutions $\theta(x,z)$ and $\xi (x,z)$:
   \[
  \{\theta(\cdot,z) , \xi(\cdot,z)\}=-1 .
\]
It follows that
 \begin{equation}
 \varphi (x,z)=  \{\varphi(\cdot,z) , \theta(\cdot,z)\} \xi (x,z)-    \{\varphi(\cdot,z) , \xi(\cdot,z)\}  \theta (x,z)
\label{eq:D6}\end{equation}
where $\{\varphi (\cdot,z), \xi (\cdot,z)\} =p(0) \xi (0,z)$. In view of Theorems~\ref{EIK} and \ref{EG},
  relation \e{eq:D6} yields the asymptotic behavior of the regular solution.

\begin{theorem}\label{EG1}
Let Assumption~\ref{PQ} be satisfied, and
let $z\in {\Bbb C} \setminus [0,\infty)$.  Then  
 \begin{equation}
 \varphi (x,z)= \frac{w(z)}{2 \sqrt{-p_{0}z}}
    a (x,z) ^{-1} (1+ o(1)), \q x\to\infty,
\label{eq:D7}\end{equation}
if $w(z) =\{\varphi(\cdot,z) , \theta(\cdot,z)\}  \neq 0$ and
 \[
 \varphi (x,z)= -  \{\varphi (\cdot,z), \xi (\cdot,z)\} a (x,z)  (1+  o(1)), \q x\to\infty,
\]
if $w(z) = 0$.
 \end{theorem}
 
 Thus $ \varphi (x,z)$ exponentially grows if $z$ is not an eigenvalue of $H$, and it exponentially decays in the opposite case.
 
 \begin{remark}\label{EGR}
   Estimates of the remainders in 
\e{eq:D4} and \e{eq:D7} are not uniform in $z$ as it approaches the half-axis $(0,\infty)$. As a consequence, we cannot put $z=\lambda\in {\Bbb R}_{+}$ in  \e{eq:D7}. Such a relation would contradict \e{eq:Ras}.
 \end{remark}

%***********************************************************
\section{Spectral results}
%* **********************************************************

 \subsection{Differential operator}

First, we define    differential operator \e{eq:A} as a self-adjoint operator.  We choose the boundary condition $f(0)=0$;
see Sect.~5.2 for other conditions.

The simplest possibility is to define $H$ via the quadratic form
\begin{equation}
h[f,f]= \int_{0}^\infty \big( p(x)|f'(x)|^2+ q(x) |f(x)|^2\big) dx.
\label{eq:Q1}\end{equation}
As is well known, this form is closed on the Sobolev space ${\sf H}^1_{0}({\Bbb R}_{+})=:{\cal D}[h]$ of functions satisfying the condition $f(0)=0$ provided
\begin{equation}
0 <p_{0}\leq p(x)\leq p_{1}<\infty \q{\rm and}\q   \sup_{x\in {\Bbb R}_+}\int_{x}^{x+1}|q(y)|dy<\infty.
\label{eq:Q}\end{equation}
Therefore $H$ may be defined as a self-adjoint operator with domain ${\cal D}(H)\subset  {\cal D}[h]$ corresponding to this form (see Chapter~10 of the book \cite{BS}). This operator satisfies the relation
 \begin{equation}
  (Hf,g)= h[f,g] 
\label{eq:Ag}\end{equation}
for all $ f\in {\cal D}(H)$ and all $g\in {\cal D}[h]$.
   
For the operator \e{eq:A}, its domain  ${\cal D}(H)$ can be described explicitly.

\begin{lemma}\label{SA}
In addition to  \e{eq:Q}, assume that the function $p(x)$ is absolutely continuous on ${\Bbb R}_{+}$.
Then $f\in {\cal D}(H)$ if and only if $f\in {\sf H}^1_{0}({\Bbb R}_{+})$, the function $p(x)f'(x)$ is absolutely continuous and
\begin{equation}
(H f)(x):=-( p(x) f' (x))' +q(x)f(x)\in   L^2({\Bbb R}_{+}).
\label{eq:Af}\end{equation}
 \end{lemma}
 
   \begin{pf}
  Suppose that $f(x)$ satisfies these conditions. Using that the function $p(x) f'(x) $ is absolutely continuous
  and    integrating by parts, we see that
\begin{equation}
  \int_{0}^\infty\big( - (p(x) f'(x))' + q(x) f(x)\Big)\,\ov{g(x)} dx= \int_{0}^\infty \big( p(x) f'(x) \ov{g'(x)} + q(x) f(x) \ov{g(x)}\big) dx
\label{eq:Aff}\end{equation}
at least for all $ g \in C_{0}^\infty({\Bbb R}_{+})$. An arbitrary $g\in {\sf H}^1_{0}({\Bbb R}_{+})$ can be approximated 
 by functions $ g_{n} \in C_{0}^\infty({\Bbb R}_{+})$ in the norm of ${\sf H}^1({\Bbb R}_{+})$.   Passing to the limit $n\to\infty$ in the equality \e{eq:Aff} for $g_{n}$ and using \e{eq:Af}, we obtain relation \e{eq:Ag} whence $f\in {\cal D}(H)$.

 Conversely, suppose that $f\in {\cal D}(H)$. It follows from \e{eq:Ag} that
 \[
| h[f,g] |\leq C \| g \|, \q \forall g\in {\cal D}[h],
\]
and hence there exists $t \in L^2({\Bbb R}_{+})$ such that
\[
 h[f,g] = \int_{0}^\infty t(x) \ov{g(x)} dx.
\]
Comparing this expression with \e{eq:Q1}, we see that
\begin{equation}
  \int_{0}^\infty  p(x)f'(x) \ov{g'(x)} dx=  \int_{0}^\infty\big( t(x) -q(x) f(x)\big) \ov{g(x)} dx.
\label{eq:Q2}\end{equation}
If $ g \in C_{0}^\infty({\Bbb R}_{+})$, then \e{eq:Q2} is the definition of the distributional derivative of the function $pf'$. This derivative equals $qf-t$ where $qf\in L_{\rm loc}^1({\Bbb R}_{+})$ because $f\in C({\Bbb R}_{+})\subset {\sf H}^1({\Bbb R}_{+})$. Therefore the function $p f'$ is absolutely continuous and
$-(p f')'+qf =t\in L^2({\Bbb R}_{+})$. 
 \end{pf}

 \subsection{Resolvent}
  Theorem~\ref{EIK}  allows us to perform spectral analysis of the operator  $H$  in a sufficiently standard way. We start with a construction of its resolvent.
 Let   $R(z) $, $\Im z\neq 0$, be an integral operator defined by the formula 
 \begin{equation}
(R( z) g)(x)=\frac{\varphi(x,z)}{w(z)}  \int_{x}^\infty \theta(y,z) g(y) dy  +\frac{\theta(x,z)}{w(z)}  \int_0^{x} \varphi(y,z) g(y) dy .
\label{eq:BRe}\end{equation}
Using that the functions $\varphi$ and $\theta$ satisfy the equation \e{eq:A1}, one easily verifies that, for example, for $g\in C_{0}^\infty ({\Bbb R}_{+})$, the function $f(x)=(R(z)g)(x)$ belongs to the domain ${\cal D}(H)$ of the operator $H$ and 
 $ (H-z)f= g$ whence $R(z)g=(H-z)^{-1}g$. It follows that $R(z)=(H-z)^{-1}$ is the resolvent   of the operator $H$. In particular, the operator   $R(z)$ is bounded.
 Relation \e{eq:BRe} means that the integral kernel of $R(z)$ equals 
\begin{equation}
R(x,y;z)=w(z) ^{-1}\varphi(x,z)\theta(y,z)\q \mbox{for} \q x\leq y \q \mbox{and}\q  R(y,x;z)= R(x,y;z).
\label{eq:B1}\end{equation}

Recall that, for $\lambda < 0$,  
the values $w(\lambda\pm i0)$ and, more generally, $\theta(x,\lambda\pm i0)$ are different, but  satisfy \e{eq:Acx}.
In particular, 
$w(\lambda + i0)=0$ if and only if $w(\lambda - i0)=0$. Nevertheless, the function $\theta(x,z)/ w(z)$ is analytic in ${\Bbb C}
\setminus [0, \infty)$ because
\begin{equation}
\frac{\theta(x,\lambda+i0)}{w(\lambda+i0)}= \frac{\theta(x,\lambda-i0)}{w(\lambda-i0)},\q \lambda<0.
\label{eq:k}\end{equation}
Indeed, consider an auxiliary function
\[
\Delta (x,\lambda)= \theta(x,\lambda+i0) w(\lambda-i0)- \theta(x,\lambda-i0) w(\lambda+i0).
\]
It satisfies equation \e{eq:A1}, $\Delta (0,\lambda)=0$ and 
\[
p(0)\Delta' (0,\lambda)=\{ \theta(\lambda+i0) ,  \theta(\lambda-i0) \}.
\]
Using asymptotics \e{eq:A22} and \e{eq:A22D}, we find that $\Delta' (0,\lambda)=0$ and hence $\Delta (x,\lambda)=0$ for all $x\geq 0$. This proves equality \e{eq:k}.

Formula  \e{eq:B1} (as well as   \e{eq:Bw1}) implies that $w(z)=0$ if and only if $z$ is an eigenvalues of the operator $H$. In particular, zeros of the function $w(z)$ are negative (it is also not excluded that $w(0)=0$). Recall also that
$w(\lambda\pm i0) \neq 0$ for $\lambda>0$.

 Let us summarize these results.

\begin{theorem}\label{R}
 Let Assumption~\ref{PQ} be satisfied. Then
 \begin{enumerate}[{\rm(i)}]
 \item
The  resolvent $R(z)=(H-z)^{-1}$ of the operator $H$ is an integral operator with kernel \e{eq:B1}.
 For all $x,y\geq 0$,  it    is an analytic function of $z\in {\Bbb C}\setminus [0,\infty)$ with simple poles at eigenvalues of the operator $H$. A point
$ z\in {\Bbb C}\setminus [0,\infty)$ is an eigenvalue of $H$ if and only if $w(z)=0$.
 \item
 For all $x,y\geq 0$,   the resolvent kernel  $R(x,y;z)$   is a continuous function of $z$ up to the cut along $[0,\infty)$ with possible exception of the point $z=0$.
 \end{enumerate}
  \end{theorem}
  
  The last assertion can also be stated as follows. Let functions $f,g \in L^2({\Bbb R}_{+})$ have compact supports. Then  the function $(R(z)f,g)$  is  continuous in $z$ up to the cut along $[0,\infty)$ with possible exception of the point $z=0$.
 This result   is known as the limiting absorption principle. It implies
 
 \begin{corollary}\label{Rc}
The spectrum of the operator $H$ on the half-axis $(0,\infty)$ is absolutely continuous.
 \end{corollary}

 %***********************************************************
\subsection{Eigenfunction expansion}

 %***********************************************************
 
Let us now consider the spectral projector $E(\lambda)$ of the operator $H$. It is also an integral operator with kernel
$E(x,y;\lambda)$  related to the resolvent kernel of $H$ by the Cauchy-Stieltjes formula
\begin{equation}
2\pi i dE(x,y;\lambda) /d\lambda=R(x,y;\lambda+i0) -R(x,y;\lambda- i 0) 
\label{eq:C}\end{equation} 
where all terms are continuous functions of $\lambda>0$.
Using now formula \e{eq:B1} for $z=\lambda\pm i0$, $\lambda>0$, $x\leq y$ and taking into account relation \e{eq:R1}, we see that  \begin{align}
  dE(x,y;\lambda) /d\lambda =& (2\pi i )^{-1} \varphi(x, \lambda)\Big( \frac{\theta(y, \lambda+i0)}{w ( \lambda+i0)}
-\frac{\theta(y, \lambda-i0)}{w ( \lambda-i0)}\Big)
\nonumber \\
=& \sqrt{ p_{0}\lambda}\:  K(\lambda)^2 \frac{\varphi (x,\lambda)\varphi (y,\lambda)}{\pi |w(\lambda \pm i0)|^2},
 \q \lambda>0.
\label{eq:C1}\end{align} 
Of course this representation extends to all $x,y\geq 0$.

The following assertion is a direct consequence of   \e{eq:C1}. It supplements the limiting absorption principle.

\begin{theorem}\label{RE}
Let   $\Sigma$ be the operator of multiplication by a function $\sigma\in L^2 ({\Bbb R}_{+})$.
Under Assumption~\ref{PQ},  the operator valued function $d \Sigma E(\lambda) \Sigma/ d\lambda$ depends continuously on $\lambda$ $($at least$)$ in the Hilbert-Schmidt norm.
   \end{theorem}

Relation    \e{eq:C1} allows us to diagonalize the operator $H$  in the same way as in the short-range case  (see, for example, \S 4.2 of the book \cite{Y}). Keeping in mind scattering theory framework, we introduce two sets of eigenfunctions $\psi_{\pm}$ and  two diagonalizations $\Psi_{\pm}$  of the absolutely continuous part $H E({\Bbb R}_{+})$ of $H$. Let us set
\begin{equation}
 \psi_{\pm}(x,\lambda)=  \frac{ \sqrt[4]{ p_{0}\lambda}\, K(\lambda)}{\sqrt{\pi}w(\lambda\mp i0)} \varphi(x,\lambda)
\label{eq:C2}\end{equation} 
and 
\begin{equation}
 (\Psi_{\pm}f)(\lambda) = \int_{0}^\infty \ov{\psi_{\pm}(x,\lambda) } f(x) dx, \q f\in L^2 ({\Bbb R}_{+})\cap
L^1 ({\Bbb R}_{+}).
\label{eq:C3}\end{equation}
It follows from \e{eq:C1} that
\begin{equation}
\|E(\Lambda)f \|^2= \int_{\Lambda} | (\Psi_{\pm}f)(\lambda)|^2 d\lambda 
\label{eq:C4}\end{equation}
 for every interval $\Lambda$ such that  
$\clos\Lambda\subset {\Bbb R}_{+}$ and hence 
\begin{equation}
 \Psi_{\pm}^* \Psi_{\pm}=E({\Bbb R}_{+} ). 
\label{eq:C5}\end{equation}
In particular, $\Psi_{\pm}$ extends to a bounded operator on $L^2 ({\Bbb R}_{+})$.

Let ${\sf A}$ be the operator of multiplication by $\lambda$ in the space $L^2 ({\Bbb R}_{+})$. Since the functions $\psi_{\pm}(x,\lambda)$ satisfy equation \e{eq:A1}, the intertwining property
\begin{equation}
 \Psi_{\pm} H= {\sf A}\Psi_{\pm}
 \label{eq:C6}\end{equation}
 holds.
 
 Relation \e{eq:C5} is naturally interpreted as the completeness of the eigenfunctions $ \psi_{\pm}(x,\lambda)$ of the operator $H$. Their orthogonality means that the adjoint operator $\Psi_{\pm}$ is isometric:
 \begin{equation}
 \Psi_{\pm} \Psi_{\pm}^*=I.
\label{eq:C7}\end{equation}
This relation can be checked exactly as in the short-range case.

Let us summarize these results.

\begin{theorem}\label{EE}
 Let Assumption~\ref{PQ} be satisfied.   Then the operators $\Psi_{\pm}: L^2 ({\Bbb R}_{+})\to L^2 ({\Bbb R}_{+})$   defined by formulas \e{eq:C2}, \e{eq:C3}   are bounded and  
 relations \e{eq:C5} -- \e{eq:C7} hold true.
 \end{theorem}
 
 \begin{corollary}\label{EEc}
 The positive spectrum of the operator $H$ covers ${\Bbb R}_{+}$. It is absolutely continuous and simple.
  \end{corollary}
 
 %***********************************************************
\subsection{Wave operators}

  %***********************************************************

Let us now consider the differential operator
\[
H_{0}=-p_{0}d^2/dx^2
\]
with the same boundary condition $f(0)=0$ in the space $L^2 ({\Bbb R}_{+})$. In this case the two operators $\Psi_{\pm}$ defined by \e{eq:C2}, \e{eq:C3} reduce to the single operator (the Fourier sine transform)
\[
 (\Psi_{0}f)(\lambda) =\frac{1}{\sqrt{\pi}  \sqrt [4]{ p_{0}\lambda }}  \int_{0}^\infty\sin \big(\sqrt{ \lambda/ p_{0}} \,x\big) f(x) dx.
\]
This operator possesses of course all properties enumerated in Theorem~\ref{EE}.

Stationary wave operators $U_{\pm}$ for the pair $H_{0}, H$ are defined by the relation
\[
 U_{\pm}=\Psi_{\pm}^* \Psi_{0}.
\]
The following result is a direct consequence of Theorem~\ref{EE}.

\begin{theorem}\label{EEW}
Under Assumption~\ref{PQ}   the wave operators $U_{\pm}$ are isometric, $ U_{\pm}^* U_{\pm}=I$, complete, i.e.,
$ U_{\pm} U_{\pm}^*= E({\Bbb R}_{+})$, and enjoy the intertwining property
$H U_{\pm} =U_{\pm} H_{0}$.
 \end{theorem}

 It follows from relations  \e{eq:C2}, \e{eq:C3} that, for all $f\in L^2 ({\Bbb R}_{+})$,
 \begin{equation}
 (\Psi_+ f)(\lambda) =S(\lambda) (\Psi_- f)(\lambda) 
 \label{eq:S}\end{equation}
 where the coefficient
 \[
S(\lambda)=\frac{w(\lambda - i0)}{w(\lambda + i0)}
\]
is known as the stationary scattering matrix.

Time-dependent wave operators $W_{\pm}$ are defined as strong limits
 \begin{equation}
 W_{\pm} = \slim_{t\to \pm\infty} e^{iHt}{\sf U}_{0} (t)
 \label{eq:WS}\end{equation}
 where ${\sf U}_{0} (t)$ is a suitable (see, for example,   \cite{BM} or \cite{Hor, Y2}) unitary regularization of $e^{-iH_{0}t}$.
 If the limits \e{eq:WS} exist, then the operators $W_{\pm}$ are automatically isometric and enjoy the intertwining property. Moreover, the operators $W_{\pm}$ are complete because the spectrum of the operator $H$ is simple according to the classical Weyl result.

Consideration of the  operators $W_{\pm}$
 is out of the scope of the present article. We  note however that various  proofs of the existence of limits in \e{eq:WS} require somewhat more stringent  conditions  on $p$ and $q$ compared to Assumption~\ref{PQ}. Under such conditions the equality $W_{\pm}=U_{\pm}$ also holds.

  %***********************************************************
\section{Miscellaneous}
%***********************************************************

     %***********************************************************
\subsection{Short-range perturbations}
%***%***********************************************************
Let us now consider a more general case where $q(x)$ is replaced by a function
 \[
q(x)+ q_{\rm sr} (x)
\]
with a short-range term $q_{\rm sr}\in L^1 ({\Bbb R}_{+})$. We suppose that $p$ and $q$ satisfy Assumption~\ref{PQ}  and define the functions $ \Omega$ and $ \theta$, etc., by formulas of Section~2 (neglecting $q_{\rm sr}$). Then for the remainder \e{eq:A3}, we have the expression
\begin{equation}
r(x,z) = (p(x) \omega(x,z))' + q_{\rm sr} (x)
\label{eq:q1}\end{equation}
instead of \e{eq:A8}; obviously, $r(\cdot,z) \in L^1 ({\Bbb R}_{+})$. Therefore
  the integral equation \e{eq:A17} with kernel \e{eq:A18} and  remainder \e{eq:q1}
 has a solution $u(x,z)$ satisfying estimate  \e{eq:A20}. Then  the Jost solution of the differential equation \e{eq:A1} with $q$ replaced by $q + q_{\rm sr}$ is again defined by equality  \e{eq:Jost}. Thus all the results obtained in the previous sections  for the particular case $q_{\rm sr}=0 $ remain true.

   %***********************************************************
\subsection{General boundary conditions}
%***%***********************************************************

Here, we briefly discuss the differential operator
\e{eq:A}
in   $L^2 ({\Bbb R}_{+})$
with a boundary condition 
\begin{equation}
f'(0)=h f(0), \q h=\bar{h}.
\label{eq:BC}\end{equation}
As before, let $\theta(x,z)$ be the Jost solution of equation \e{eq:A1}
 constructed in Theorem~\ref{EIK}. Instead of \e{eq:B}, the regular solution $\varphi (x,z)$  of this equation will now be distinguished by the conditions
 \begin{equation}
\varphi(0,z)=1, \q \varphi'(0,z) =h.
\label{eq:BC1}\end{equation} 
This is again an analytic function of $z\in{\Bbb C}$. Formula \e{eq:R1} remains true where according to \e{eq:BC1}  the Wronskian
 \begin{equation}
w(z)=\{ \varphi(\cdot,z), \theta(\cdot,z) \} = p(0) (h \theta(0,z)-\theta'(0,z)).
\label{eq:BC2}\end{equation}

The resolvent kernel is still given by the relation \e{eq:B1} which yields representation \e{eq:C1} with the functions $ \varphi(x, \lambda)$ and $w(z)$  defined by  \e{eq:BC1}  and  \e{eq:BC2}. As before,  the functions $\psi_{\pm} (x,\lambda)$ are given by formula \e{eq:C2} and the operators $\Psi_{\pm}$ -- by formula \e{eq:C3}. Then Theorems~\ref{EE} and \ref{EEW} remain true for the operator $H$ corresponding to the boundary condition \e{eq:BC}.

%***********************************************************
\subsection{Problem on  the whole line}
%***%***

 Consider now the operator \e{eq:A} in the space $L^2({\Bbb R})$. We here follow closely the scheme described for short-range potentials, for example, in \S 5.1 of \cite{Y}.
 
  Suppose that the conditions on $p(x)$ and $q(x)$ are imposed for all $x\in {\Bbb R}$; in particular, 
   the limits in \e{eq:A2a} are taken for $|x|\to\infty$. In addition to the Jost solution $\theta(x,z)=:\theta_1(x,z)$ of equation \e{eq:A1}  built in Theorem~\ref{EIK}, we distinguish a solution $\theta_2(x,z)$ by its asymptotics as $x\to-\infty$. For the function $\Omega(x,z)$ defined by 
\e{eq:A11}, we set
$
a_2(x,z) =e^{\Omega(x,z)}
$. Then the construction of Theorem~\ref{EIK} leads to the solution  $\theta_2(x,z)$ of equation \e{eq:A1}
 with asymptotics     $\theta_2(x,z)\sim a_2(x,z)$ as $x\to -\infty$ (cf. \e{eq:A22}).
 We also introduce   the Wronskian of the solutions $\theta_{1}$ and $\theta_{2}$:
\[
 w(z)=\{\theta_{2}(\cdot,z) , \theta_1(\cdot,z)\}.
\]
 The construction of the resolvent $R(z)=(H-z)^{-1}$, $z\in \Pi$,  is similar  to Subs.~4.2.   Since $\theta_{1}\in L^2$ as $x\to\infty$ and  $\theta_2\in L^2$ as $x\to -\infty$, the  resolvent   kernel equals (cf. \e{eq:B1})
\begin{equation}
R(x,y;z)=w(z) ^{-1}\theta_{2}(x,z)\theta_{1}(y,z)\q \mbox{for} \q x\leq y \q \mbox{and} \q  R(y,x;z)= R(x,y;z).
\label{eq:B1w}\end{equation}

Suppose now that $z=\lambda\pm i0$ where $\lambda>0$. Calculating the Wronskians for $x\to \infty$ (if $j=1$) or for $x\to -\infty$ (if $j=2$), we find that
  \[
  \{\theta_{j}(\lambda+ i0), \theta_{j}(\lambda- i0)\}=(-1)^{j-1}2 i \sqrt{p_{0}\lambda}, \q j=1,2.
  \]
  Thus, these solutions are linearly independent for all $\lambda>0$, and we have
  \begin{align}
   \theta_1 (x, \lambda+ i0)&= (2 i)^{-1} (p_{0}\lambda)^{-1/2} K_{1}(\lambda)^{-2} \big(\ov{{\bf w}(\lambda)}   \theta_2 (x, \lambda+ i0)
-  w (\lambda+i0)   \theta_2 (x, \lambda- i0) \big)
  \label{eq:px1}
  \\
  \theta_2 (x, \lambda+ i0)&= (2 i)^{-1} (p_{0}\lambda)^{-1/2}K_2(\lambda)^{-2}\big({\bf w}(\lambda)   \theta_1 (x, \lambda+ i0)
-  w (\lambda+i0)   \theta_1 (x, \lambda- i0) \big)
   \label{eq:px2}\end{align}
 where
  \[
  {\bf w}(\lambda) =\{\theta_2 (\cdot, \lambda+ i0), \theta_1 (\cdot, \lambda- i0)\},
   \]
     $  K_{1}(\lambda)= K (\lambda)$ is defined by \e{eq:K} and
      \[
  K_{2}(\lambda)= \exp\Big(-\int_{-\infty}^ {0}  \sqrt{\Big( \frac{ q(y)-\lambda}{p(y)}\Big)_{+}}dy\Big).
 \]
 
   Note the identity
       \begin{equation}
  |w (\lambda\pm i0) |^2= 4 p_{0} \lambda K_1(\lambda)^{2}K_2(\lambda)^{2}+ |{\bf w}(\lambda)|^2.
   \label{eq:ph2}\end{equation}
   For its proof, we
        substitute   expression \e{eq:px1} for $ \theta_1 (x, \lambda+ i0)$  and 
        $ \ov{\theta_1 (x, \lambda+ i0) }= \theta_1 (x, \lambda- i0)$ into the right-hand side of \e{eq:px2} and observe that the coefficient in front of  $ \theta_2 (x, \lambda+ i0)$ should be equal to $1$.

We need an analogue of representation \e{eq:C1}.

\begin{lemma}\label{Exw}
For all $x,y\in{\Bbb R}$ and $\lambda>0$, we have the representation
  \begin{equation}
 dE(x,y;\lambda) /d\lambda=  \frac{ \sqrt{ p_{0}\lambda} }{\pi |w(\lambda \pm i0)|^2}
\Big( K_1(\lambda)^{2} \theta_{1} (x,\lambda+i0) \theta_{1} (y,\lambda-i0) + K_2(\lambda)^{2}\theta_2 (x,\lambda+i0) \theta_2 (y,\lambda-i0) \Big).
\label{eq:C1w}\end{equation} 
\end{lemma}

\begin{pf}
It follows from  the Cauchy-Stieltjes formula \e{eq:C} and the representation \e{eq:B1w}  that 
  \begin{multline}
 dE(x,y;\lambda) /d\lambda=  \frac{ 1} {2\pi i |w(\lambda \pm i0)|^2}
\Big( w(\lambda - i0)\theta_{2} (x,\lambda+i0) \theta_{1} (y,\lambda+i0)
\\
-w(\lambda + i0)\theta_{2} (x,\lambda-i0) \theta_{1} (y,\lambda-i0) \Big), \q x\leq y.
\label{eq:Cxw}\end{multline} 
Let us show that the right-hand sides of \e{eq:C1w} and \e{eq:Cxw} coincide. Replacing $\theta_{1} (x,\lambda+i0)$ in \e{eq:C1w} by
its expression \e{eq:px1}, we see that it suffices to check the identity
  \begin{multline}
\big({\bf w}(\lambda) \theta_2(x,\lambda+i0) - w(\lambda+i0) \theta_2 (x ,\lambda-i0) \big) \theta_1 (y ,\lambda-i0)+ 2i \sqrt{ p_{0}\lambda} \, K_2(\lambda)^{2}\theta_2 (x,\lambda+i0) \theta_2 (y,\lambda-i0) 
 \\
 = w(\lambda - i0)\theta_{2} (x,\lambda+i0) \theta_{1} (y,\lambda+i0) -
w(\lambda + i0)\theta_{2} (x,\lambda-i0) \theta_{1} (y,\lambda-i0).
 \label{eq:CXW}\end{multline}
    The coefficients in front of $\theta_2 (x, \lambda+i0)$ in the left and right sides of \e{eq:CXW} coincide by virtue of 
    \e{eq:px2}. The terms containing $\theta_2 (x, \lambda-i0)$ in the left and right sides of \e{eq:CXW} are the same.
Of course representation  \e{eq:C1w} extends to $x \geq y$ since its right-hand side becomes complex conjugated if $x$ and $y$ are interchanged.
\end{pf}
 
 \begin{remark}\label{ERw}
 The left-hand side of \e{eq:C1w} is real and symmetric in $(x,y)$. Therefore in addition to \e{eq:C1w}, we have the representation
  \[
 dE(x,y;\lambda) /d\lambda=  \frac{ \sqrt{ p_{0}\lambda} }{\pi |w(\lambda \pm i0)|^2}
\Big( K_1(\lambda)^{2} \theta_{1} (x,\lambda-i0) \theta_{1} (y,\lambda+i0) + K_2(\lambda)^{2} \theta_2 (x,\lambda-i0) \theta_2 (y,\lambda+i0) \Big).
\]
\end{remark}
 
Instead of \e{eq:C2}, we now define eigenfunctions of the continuous spectrum of the operator $H$ by the relation
\begin{equation}
 \psi_j(x,\lambda)=  \frac{ \sqrt[4]{ p_{0}\lambda}\, K_{j} (\lambda) }{i \sqrt{\pi}w(\lambda + i0)} \theta_{j}(x,\lambda+ i0), \q j=1,2.
\label{eq:C2w}\end{equation} 
Then (cf.  \e{eq:C3}) we introduce 
  the mappings $\Psi_{\pm} : L^2 ({\Bbb R})\to L^2 ({\Bbb R}_{+}; {\Bbb C}^2)$ by   formulas
\begin{equation}
\begin{split}
 (\Psi_{+}f)(\lambda)& =\Big (\int_{-\infty}^\infty\psi_2(x,\lambda) f(x) dx ,
 \int_{-\infty}^\infty\psi_1 (x,\lambda) f(x) dx\Big)^\top
 \\
 (\Psi_{-}f)(\lambda) &=\Big (\int_{-\infty}^\infty \ov{\psi_1(x,\lambda) }f(x) dx ,
 \int_{-\infty}^\infty \ov{\psi_2(x,\lambda) }f(x) dx\Big)^\top.
\end{split}
\label{eq:C3w}\end{equation}

It follows from \e{eq:C1w} that relation \e{eq:C4} holds. Therefore,
similarly to the proof of Theorem~\ref{EE}, one obtains the following result. Note that the multiplication operator $\sf A$ acts now in the space $L^2 ({\Bbb R}_{+}; {\Bbb C}^2)$.

\begin{theorem}\label{EEw}
 Let Assumption~\ref{PQ} be satisfied for all $x\in{\Bbb R}$.
 Then the operators $\Psi_{\pm} $ are bounded and satisfy
 relations \e{eq:C5} -- \e{eq:C7}.
 \end{theorem}
 
 \begin{corollary}\label{EEcw}
 The positive spectrum of the operator $H$ covers ${\Bbb R}_{+}$. It is absolutely continuous and  has multiplicity two.
  \end{corollary}

   In terms of functions  \e{eq:C2w}, relations \e{eq:px1} and \e{eq:px2} can equivalently be rewritten as
    \begin{equation}
   \begin{pmatrix}
\psi_2 (x, \lambda ) \\ \psi_1 (x, \lambda )
    \end{pmatrix}
    = S  (\lambda)
  \begin{pmatrix}
\ov{\psi_1 (x, \lambda ) }\\  \ov{ \psi_2 (x, \lambda )}
    \end{pmatrix}
\label{eq:SS}\end{equation}
where  
    \begin{equation}
S  (\lambda)=   w(\lambda+ i0)^{-1}\begin{pmatrix}
 i \gamma (\lambda)& {\bf w}(\lambda) \\
\ov{{\bf w}(\lambda) }& i \gamma (\lambda)      \end{pmatrix} \q\mbox{and}\q \gamma(\lambda)= 2\sqrt{p_{0}\lambda} K_{1}(\lambda)
K_2(\lambda).
\label{eq:SS1}\end{equation}
 According to \e{eq:ph2} the $2\times 2$ matrix $S  (\lambda)$   is unitary.
 It is  known as the scattering matrix for the problem on the whole line. 
It follows from \e{eq:SS}   that, for all $f\in L^2 ({\Bbb R})$, the identity 
  \e{eq:S}  holds with the operators $\Psi_{\pm}$ defined by \e{eq:C3w}  and the matrix  \e{eq:SS1}.

Finally, we note that formulas \e{eq:px1}, \e{eq:px2} yield asymptotics of the eigenfunctions \e{eq:C2w} as $x\to\pm\infty$.

\end{document}